\documentstyle[11pt]{article}
\newcommand{\qed}{\quad \rule[-.5mm]{1.9mm}{3.2mm}}
\newcommand{\equsep}{\rule{0mm}{6mm}}

\newcommand{\newsection}[1]{\section{#1} \setcounter{equation}{0}}

\newtheorem{thm}{Theorem}[section]
\newtheorem{prop}[thm]{Proposition}
\newtheorem{de}[thm]{Definition}
\newtheorem{lem}[thm]{Lemma}

\begin{document}

 \setlength{\baselineskip}{1.25\baselineskip}

\begin{center}
{\large\bf Restricted Routing and Wide Diameter of the Cycle
		Prefix Network}
\end{center}

\begin{center}
\parbox{15cm}{\begin{center}
\baselineskip=1.1\baselineskip 
    William Y. C. Chen, $\;$  Vance Faber $\;$  and $\;$ 
		Emanuel Knill

     \vspace{3mm}
      {\em C-3, Mail Stop B265\\
        Los Alamos National Laboratory\\
        Los Alamos, New Mexico 87545}

        \end{center}
     }
\end{center}

\vspace{3mm}

\begin{abstract}
The cycle prefix network is a Cayley coset digraph based
on sequences over an alphabet which has been proposed
as a vertex symmetric communication network.
This network
has been shown to have many remarkable communication
properties such as a large number of vertices for a given degree
and diameter, simple shortest path routing, Hamiltonicity,
optimal connectivity, and others.
These considerations for designing symmetric
and directed interconnection networks are well justified in practice
and have been widely recognized in the research community.
Among
the important properties of a good network,
efficient routing is probably one of the most important. In this 
paper, we further study routing schemes in the cycle prefix
network. We confirm an observation first made from computer
experiments regarding the diameter change when certain
links are removed in the original network, and we completely
determine the wide diameter of the network.
The wide diameter of a network is now perceived to be even more important
than the diameter.
We show by construction that the wide diameter of the cycle
prefix network is very close to the ordinary diameter.
This means that routing in parallel in this network
costs little extra time compared to ordinary 
single path routing. 
\end{abstract}

\vspace{2cm}

\begin{center}

   Suggested Running Title:

   Cycle Prefix Network

        .................................................


\end{center}

\newpage

\newsection{Introduction}

The cycle prefix network is a vertex symmetric, directed
graph which has recently been proposed for use as a communication
network [\ref{FMC93}]. It has been shown that the cycle prefix
network has many
remarkable communication properties such as a large number of
vertices for a given degree and diameter, simple shortest path routing,
Hamiltonicity, optimal connectivity, and others 
[\ref{FMC93},\ref{CFK93},\ref{JR92}].
These considerations
for designing symmetric and directed networks are well justified in
practice and have been widely recognized in the research community.
In the search for highly efficient network models and in the
study of communication algorithms on these networks, Cayley
graph techniques have been used successfully in discovering new models and
in analyzing network efficiencies.
Another interesting idea
in network design which is used for the construction of many networks
is to represent nodes as
certain sequences over an alphabet with
links represented by suitable operations on sequences. The
hypercube, de Bruijn, Kautz, star and pancake 
networks can all be constructed in this fashion.
In the case of cycle prefix digraphs, both the
idea of Cayley graphs (Cayley coset digraphs, to be precise)
and that of sequences over an alphabet can be used
as the underlying representation, and each has
its own advantage. The former idea 
was the point of departure for the discovery of the
cycle prefix network, motivated by
the fundamental theorem of Sabidussi [\ref{Sab69}] which shows that
all vertex symmetric digraphs are Cayley coset digraphs.
The sequence representation of a cycle prefix network
is more useful for studying its properties and for implementing
it in practice.
For this reason,
we shall utilize the sequence representation of the
cycle prefix digraph throughout the paper. 

In the performance evaluation of networks, the diameter and routing
efficiency
are among the most critical concerns.
In this paper we shall
further study these issues for the cycle prefix network. 
The paper has two objectives.
The first objective is to confirm
an observation first made on the basis of computer experiments concerning
the diameter change when certain links are removed. 
By proving a reachability theorem we show that
in the case of cycle prefix digraphs, the resulting networks often
possess better degree-diameter properties than the original
one. Using the method of Conway and Guy [\ref{CG}], 
one may construct large symmetric
networks with small degree and diameter.
In a recent
paper [\ref{CF92}], Comellas and Fiol describe a routing scheme
with implies an upper bound on the diameter of the link deleted
cycle prefix digraphs. However, they left open the question
of whether the bound is exact. We settle this question by 
exhibiting vertices that achieve the diameter bound. 

The second objective of this paper is concerned with the
recently introduced notion of the wide diameters of networks.
Mathematically,
the notion of the wide diameter of a graph naturally 
stems from the classical theorem of Menger relating
connectivity to disjoint paths. However, such a notion
has not been studied in graph theory until very
recently when it became relevant from an engineering
point of view.
D. F. Hsu gives a vivid account on the
background of wide diameters [\ref{Hsu**}]:

``The concept of the wide diameter of a graph 
$G$ arises naturally from the study of routing,
reliability, randomized routing, fault tolerance,
and other communication protocols (such as the byzantine algorithm)
in parallel architecture and distributed computer 
networks. By considering both the width and the length
of a container, we are able to give a global and systematic 
treatment on the interconnection network for
various distributed systems.''

``Although the concept of a container and the notion of
wide diameter have been discussed and used in practical
applications, the graph theory questions suggested
have not, at least until recently, been studied
as extensively as the questions in the hardware and software
design, development, and implementation of distributed
computing systems.''

In the case of the cycle prefix networks, we completely
determine its wide diameter through an explicit construction.
It turns out that the wide diameter is very close to the 
ordinary diameter of the network. In other words, routing
in parallel in the cycle prefix network costs 
little extra time compared to ordinary single path routing.
This property undoubtedly increases the usefulness of the cycle
network.

\newsection{The Cycle Prefix Digraphs $\Gamma_\Delta(D)$
		and $\Gamma_\Delta(D, -r)$}

\newcommand{\cp}{$\Gamma_\Delta(D)$}
\newcommand{\cpo}{$\Gamma_\Delta(D,-1)$}
\newcommand{\cpr}{$\Gamma_\Delta(D, -r)$}

The cycle prefix digraph $\Gamma_\Delta(D)$ $(\Delta\geq D)$
is defined as a digraph whose vertex set consists of sequences
$x_1x_2\cdots x_D$ over an alphabet $\{ 1, 2, \ldots, \Delta+1\}$,
where $x_1, x_2, \ldots, x_D$ are distinct. Such sequences are
called partial permutations or $D$-permutations. The adjacency
relations for a vertex $x_1x_2\cdots x_D$ are described as
follows:
\[ x_1x_2\cdots x_D \Rightarrow \left\{ 
		\begin{array}{ll}
		  x_k x_1\cdots x_{k-1} x_{k+1} \cdots x_d, \quad
		  	& \mbox{for $2\leq k \leq D$}, \\
		\equsep
		  y x_1 x_2 \cdots x_{D-1}, &
		  	\mbox{for $y \not= x_1, x_2, \ldots, x_D$}.
		\end{array}
		\right. \]
We say that the vertex $x_kx_1\cdots x_{k-1} x_{k+1} \cdots x_D$
is obtained from $x_1x_2\cdots x_D$ via a rotation on the
prefix $x_1x_2\cdots x_k$, denoted,
\[ x_k x_1\cdots x_{k-1} x_{k+1} \cdots x_D
	= R_k (x_1x_2\cdots x_D)\; .\]
In particular, the operation $R_D$ is called a full
rotation, and $x_D x_1 x_2\cdots x_{D-1}$ is called
a full rotation of $x_1x_2\cdots x_D$. 
The rotations $R_k$ for $k< D$
are called partial rotations. Similarly, we say
that the sequence $yx_1x_2\cdots x_{D-1}$ is obtained
from $x_1x_2\cdots x_D$ via a shift, denoted
\[ yx_1 x_2\cdots x_{D-1} = S_y ( x_1x_2\cdots x_D).\] 
Since the term ``adjacent'' is somewhat ambiguous
for a directed graph, when $(u, v)$ is an arc in
a digraph we shall say that $u$ is {\em adjacent to}
$v$ while $v$ is {\em next to } to $u$. The term
``adjacent from'' is used by some authors to distinguish
it from ``adjacent to''. From the above sequence definition
of $\Gamma_\Delta(D)$, it is easily seen to be vertex
symmetric, a fact that immediately follows from the
Cayley coset digraph definition $\Gamma_\Delta(D)$
[\ref{FMC93}]. 

We notice that some authors prefer the sequence
shift in the direction from right to left, like the
shift for de Bruijn digraphs:
$(x_1, x_2, \ldots, x_n) \Rightarrow (x_2, \ldots, x_n, z)$.
For the sake of consistency,
we shall follow the notation in [\ref{FMC93}], and
continue the left-to-right shift, which reflects
the rotations on the prefixes. A left-handed notation is
however adopted by Comellas and Fiol [\ref{CF92}].

In the original study of cycle prefix digraphs $\Gamma_\Delta(D)$,
Faber and Moore observed from computational experiments that
if one rules out the double arcs in $\Gamma_\Delta(D)$,
then the diameter of the resulting digraph, denoted, \cpo, 
increases only by one. They then considered a new
construction based on the cycle prefix digraph \cp. Suppose \cpr\ 
is the digraph obtained from \cp\ by deleting the 
arcs represented by the partial rotations
$R_2, R_3, \ldots, R_{r+1}$. Formally speaking,
\cpr\ has the same vertex set as \cp, and the
adjacency relations for a vertex $x_1x_2\cdots x_D$
are described by
\[ x_1x_2\cdots x_D \Rightarrow \left\{ 
		\begin{array}{ll}
		  x_k x_1\cdots x_{k-1} x_{k+1} \cdots x_d, \quad
		  	& \mbox{for $r+2\leq k \leq D$}, \\
		\equsep
		  y x_1 x_2 \cdots x_{D-1}, &
		  	\mbox{for $y \not= x_1, x_2, \ldots, x_D$}.
		\end{array}
		\right. \]
Note that the degree of \cpr\ decreases by $r$ compared with
\cp. There is an intuitive reason to surmise that the 
diameter of \cpr\ would increase by the same scale.
Recently, Comellas and Fiol [\ref{CF92}] have shown that in most cases
the diameter increase of \cpr\ is in fact bounded by $r$. However,
they did not solve the problem of whether this bound is
exact or not. We will fill this gap
by proving the exactness of the diameter bound. Moreover,
we shall study the reachability property of the digraph
\cpr\ (this property for \cp\ has been studied 
in [\ref{CF92}]). A vertex symmetric digraph with
the reachability property is of great use in constructing
new classes of vertex symmetric digraphs with small 
degree and diameter, as proposed by Conway and Guy [\ref{CG}].
The details are presented in the next section.

\newsection{Restricted Routing for \cpr}

For the sake of easier presentation, we shall start with the
reachability of the digraph \cpr. It is easy to see
that \cpr\ is vertex symmetric and that one may choose
the {\em standard origin} $12\cdots D$ while considering
routing for any two vertices. For simplicity the vertex
$X$ is always referred to as $x_1x_2\cdots x_D$.

\begin{de}[$k$-Reachable digraphs]
A digraph is said to be $k$-reachable if
for any two vertices $u$ and $v$, which are not necessarily distinct,
there exists a path (with repeated vertices and
arcs allowed) from $u$ to $v$  of length $k$.
\end{de}

Comellas and Fiol [\ref{CF92}] have shown that the digraph \cp\  
is $D$-reachable for $\Delta \geq D \geq 3$. Here we will
present a stronger result for \cpr.  

\begin{thm}
Suppose $r\geq 0$ and $\Delta\geq D \geq 2r+3$. Then
the vertex symmetric digraph \cpr\ is $(D+r)$-reachable.
\label{thm-reach}
\end{thm}

The major reason for the above theorem lies in the
following observation about  {\em dead angles}. 
Given a vertex $X=x_1x_2\cdots x_D$, the prefix
$x_1x_2\cdots x_{r+1}$ is called the dead angle of
$X$ in \cpr. We say that a letter $z$
is in the dead angle of $X$ if $z=x_i$ for some $1\leq i \leq r+1$.

\begin{lem}[Dead Angle Principle]
Let $X$ be a vertex in \cpr. Then there exists a vertex
$Y$ next to $X$ such that $Y$ begins with a letter $z$,
if and only if $z$ is not in the dead angle of $P$.  
\end{lem}

The proof of the above lemma is straightforward.
It is based on a property of \cp\ regarding
how the rotation and shift operation work
together to complement each other:
Suppose $z$ is not in the dead angle of $X$.
If $z$ is indeed in $X$, then a rotation operation
on $X$ may put $z$ back to the beginning of the
sequence; otherwise, a shift operation can achieve
the same goal with ease. For the above reason,
one sees that the two operations are coherent with
each other, although the look rather unrelated. Moreover,
suppose $Y$ is next to $X$ in \cpr, then
$Y$ is determined by its first element. 
We now give the proof of Theorem \ref{thm-reach}.

{\em Proof.} Let $I=12\cdots D$ be the standard origin
and $X=x_1x_2\cdots x_D$ be the destination. It suffices
to show that there is a directed path of length $D+r$ from $I$
to $X$. We first consider the case when $x_D\not=1$.
Let $A$ be the set $\{x_{1}, x_{2}, \ldots, x_{D-r-1}\}$.
Since there are $r+1$ elements in the dead angle of $X$,
and $D-(r+1) > r+1$, there exists an element $y_1$
in $A$ that is not in the dead angle of $I$.
By the dead angle principle, $I$ is adjacent to a
vertex $P_1$ with prefix $y_11 \cdots (r+1)$. 
Let $A=A \backslash \{y_1\}$. Considering the
dead angle of $P_1$, the same condition on $D$
and $r$ ensures that that
 there exists an element $y_{2}$ in $A_1$
that is not in the dead angle of $P_1$ (implying that
$y_1$ and $y_2$ are distinct). Hence $P_1$ is adjacent
to a vertex $P_2$ with prefix $y_{2} y_1 1 \cdots r$. 
Repeating the above procedure, one may reach a
vertex $P_r$ such that 
$P_r$ has prefix $y_r\cdots y_2y_1 \,1$ and $y_1, y_2, 
\ldots,y_r$ come from $A$.  It has already taken $r$
steps to get to $P_r$ from $I$.  Since
$x_D\not= 1$, we may construct a path
of length $D$ from $P_r$ to $X$, and display it by
showing the prefixes:
\[ y_r \cdots y_2 y_1 1 \quad \Rightarrow \quad
	x_D y_r\cdots y_2 y_1 \quad \Rightarrow \quad
	x_{D-1} x_D y_r\cdots y_2 y_1 \quad \Rightarrow \quad \cdots 
	\quad \]
\[ \Rightarrow  \quad x_{D-r} \cdots x_{D-1} x_D y_r\cdots y_2 y_1 \quad
	\Rightarrow \quad 
   x_{D-r-1} x_{D-r} \cdots x_{D-1} x_D  \quad \Rightarrow \quad
   	\cdots \quad \] 
\[\Rightarrow \quad x_1 x_2\cdots x_D\;. \]

We next consider the case when $x_D=1$. Let $A=\{x_1, x_2, \ldots,
x_{D-r-2}\}$ and $B=\{ 2, 3, \ldots, r+1\}$. Since $D-r-2 >r$,
there exists an element $y_1$ such that $y_1 \in A$ but $y_1\not\in B$.
Note that $1\not\in A, B$. Hence by the dead angle principle,
$I$ is adjacent to a vertex $P_1$ with prefix $y_1 1 \,2\, \cdots (r+1)$.
Let $A_1=A\backslash \{y_1\}$, $B_1=B\backslash \{r+1\}$. The
same condition on $D$ and $r$ ensures that there exists $y_2\in A_1$,
but $y_2\not\in B_1$. It follows that $P_1$ is adjacent to 
a vertex $P_2$ with prefix $y_2y_1 \, 1\, 2\, \cdots \, r$.
Repeating this procedure, one ends up with a vertex $P_r$
having prefix $y_ry_{r-1}\cdots y_1 1$, where $y_i\in A$. 
We continue with the following path of length $D-r$ starting
from $P_r$ (with only prefixes shown):
\begin{eqnarray*}
   y_r y_{r-1} \cdots y_1 1   & \Rightarrow  &  
   			x_{D-r-1}   y_r y_{r-1} \cdots y_1 1  \\
     & \Rightarrow  &  x_{D-1}  x_{D-r-1}   y_r y_{r-1} \cdots y_1 1  \\
     & \Rightarrow  & x_{D-2} x_{D-1}  x_{D-r-1}   y_r y_{r-1} \cdots y_1 1 
         \\   
    & \cdots & \\
 & \Rightarrow  &  x_{D-r} \cdots
 		x_{D-1}  x_{D-r-1}   y_r y_{r-1} \cdots y_1 1  \\
 & \Rightarrow  &  x_{D-r-1} x_{D-r} \cdots
 		 x_{D-1}   y_r y_{r-1} \cdots y_1 1  
\end{eqnarray*}
The last vertex is labeled by $P_{2r+2}$ according to the
length. Note that $y_i \in A$, we claim that there is
a path from $P_{2r+2}$ to $X$ of the following form:
\[ P_{2r+2} \quad \Rightarrow \quad 
	x_{D-r-2} x_{D-r-1} \cdots x_{D-1} \cdots 1 
      \quad \Rightarrow \quad x_{D-r-3} \cdots x_{D-1} \cdots 1 
      \quad \Rightarrow \quad \cdots \]
\[ \quad \Rightarrow
      	\quad x_1x_2 \cdots x_{D-1} 1 =X\, , \]
because $y_i \in A$ at each step it is impossible to bump
1 out of the sequence so that the last vertex has to
be $X$. Summing up all the segment, we get a path of 
length $D+r$. \qed

Specializing  the above theorem for $r=0$,
it follows the reachability of \cp\ 
first observed in [\ref{CF92}]. Moreover, 
using the method of Conway and Guy [\ref{CG}], one may construct
large symmetric digraphs with small degree and diameter
based on \cpr. 
Since the digraph \cpr\ in some cases has more vertices
than $\Gamma_{\Delta-r}(D+r)$, one may use the above theorem
in constructing new symmetric digraphs. However, we will
not discuss this aspect here.

The rest of this section is concerned with the diameter of
\cpr. It is shown in [\ref{CG}] that the diameter of
\cpr\ does not exceed $D+r$ for $\Delta\geq D \geq 2r+2$.
This upper bound is achieved by the following construction that is
a much simpler version than the construction for the
reachability of \cpr. Note that the reachability result
requires a slightly stronger condition on the parameters
of \cpr. Let's give an outline. 
Let $I=12\cdots D$ be the standard origin, and $X=x_1x_2\cdots x_D$
be any vertex in \cpr. For the case $x_D\not=1$, one may
first try to reach from $I$ 
a vertex with prefix $y_ry_{r-1}\cdots y_1 1$, where $y_i \not=
x_{D-1}, x_{D-2}, \ldots, x_{D-r}$. Then one may continue with
vertices having prefixes $x_D, x_{D-1}$, etc. For the case
$x_D=1$,  one may get to a vertex with prefix
$y_r y_{r-1}\cdots y_1 1$ such that $y_i \in \{
x_1, x_2, \ldots, x_{D-r-2}\}$.  Then one may 
get to $X$ via vertices with prefixes
$x_{D-1}$, $x_{D-2}x_{D-1}$, etc. The last element $x_D=1$ will
eventually takes care of itself for the reason given in the
proof of the reachability theorem. 

It is harder to show that the above diameter bound 
is exact. For this purpose, we find a class of vertices
that achieve the bound. A vertex $X$ in \cpr\
is called a remote vertex if $x_{D-1}=1$ and $x_D=D$,
and there exists $x_i>D$ for some $1\leq i \leq D-2$. 
We shall use the common notation $d(X, Y)$ to denote
the distance from $X$ to $Y$ in a digraph.
Then
we have the following theorem:

\begin{thm}
Let $\Delta\geq D \geq 2r+2$, and $X$ be a remote vertex
in \cpr, then the distance from the standard origin
$I=12\cdots D$ to $X$ equals $D+r$.
\end{thm} 

{\em Proof.} The diameter upper bound is already established,
so it suffices to show that $d(I,X)\geq D+r$. 
Since $X=x_1x_2\cdots x_{D-2}1D$ and there exists
$x_i>D$ for some $i$, to reach $X$ from $I$ requires at least one
shift operation. Thus, element $D$ in $I$ cannot remain in the
last position during the process to reach $X$ from $I$. Since
$D$ is in the destination vertex $X$, it is either moved back to
the beginning position at some point, or is removed out of the
sequence by a shift operation and then put back to the
beginning by another shift operation. Let $I\Rightarrow
P_1 \Rightarrow P_2 \Rightarrow \cdots \Rightarrow P_m =X$
be a shortest path from $I$ to $X$. Since either a 
shift or rotation operation on a vertex, say $Y=y_1y_2\cdots y_D$,
moves the elements in the dead angle to the positions on
the right hand side, the vertex $P_{r+1}$ must have the
prefix $z_{r+1} \cdots z_2 z_1 1$. If $D$ does not appear
in $z_{r+1} \cdots z_2 z_1$, then it will take at least $D$
steps to reach $X$ from $P_{r+1}$ as far as the last element
of $X$ is concerned. This contradicts the upper bound $D+r$
on the diameter of \cpr. We now assume that $z_i=D$.
It follows that $P_{r+i+1}$ has a prefix of the form
$w_i  \cdots w_2  w_1 \, z_{r+1} \cdots z_{i+1} D$. Consider
the following two cases:

Case 1. The element 1 appears in $w_i \cdots w_2 w_1$.
        If $D$ is shifted out of a vertex after the vertex $P_{r+i+1}$,
then it will take at least $D$ steps to put $D$ to the last position
of $P_m$, a contradiction. Thus, $D$ has to remain in the vertices
$P_{r+i+1}$, $P_{r+i+2}$, $\ldots$, $P_m$. Moreover,
$D$ will never be put back to the beginning of a vertex by a
rotation because after that rotation one needs at least
$D-1$ steps to move $D$ to the last position of $P_m$, which
is also impossible. Hence, in the path from $P_{r+i+1}$ to $P_m$,
$1$ has to remain in all the vertices on this path segment,
and $1$ is always to the left of $D$. We define $\delta(Y)$
to be the number of elements between $D$ and 1. 
Let $f$ be the number of operations used to reach
$P_m$ from $P_{r+i+1}$ that move the element $D$ to its
right, and $g$ be other operations used in the same path.
For a rotation or a shift operation $T$ on $P_j$ $(r+i+1 \leq j \leq m-1)$,
if $T$ moves $D$ to its right, then $T$ leaves the value of $\delta(P_j)$
unchanged; otherwise $T$ may reduce the value of 
$\delta(P_j)$ at most by one. It follows that
\[m-r-i-1 = f+g   \geq (D-r-2) + (r-i+1) = D-i-1\, .\]
Hence $m \geq  D+r$.

Case 2. The element $1$ does not appear in $w_i \cdots w_2 w_1$. 
 	As in Case 1, the element $D$ has to remain in the vertices
	on the path from $P_{r+i+1}$ to $P_m$, and $D$
	is never moved back to the beginning of any 
	vertex on the path. It is clear that at some step,
	$1$ has to be put back to the beginning of a vertex
	on the aforementioned path either by a rotation or a shift
	operation. Suppose this happens to $P_j= 1 \cdots$ $(j\geq
	r+i+2$). Since 1 is never moved back to the beginning
	of a vertex, in $P_j$ there are at least $r+1$ elements
	between $1$ and $D$. Thus, we need at least $r+1$ steps
	to move $1$ next to $D$. It follows that
\[ m \geq (r+i+1) + (D-r-2)+(r+1)  = D+r+i \geq D+r \; .\]
This completes the proof. \qed

We remark that when $r\geq 1$ the digraph \cpr\ 
does not have the unique shortest path property like
\cp. For example,
in $\Gamma_4(4,-1)$ there are two shortest paths from 1234 to
5214, shown below:
\[ \begin{array}{ccccccccccc}
   1234 & \Rightarrow & 4123 & \Rightarrow & 5412 &
   	\Rightarrow &  1542  & \Rightarrow  & 2154  & \Rightarrow
			& 5214, \\
\equsep
   1234 &  \Rightarrow & 5123 & \Rightarrow & 4512 & \Rightarrow
          &  1452 & \Rightarrow & 2145 & \Rightarrow & 5214. 
\end{array} \]

\newsection{The Wide Diameter of \cp}

Connectivity considerations of a network was primarily motivated
by its fault tolerance capabilities, while the diameter is 
a measurement of routing efficiency along a single path. 
Interestingly, the recent notion of wide diameter is
a kind of unification of both the diameter and the
connectivity due to the classical theorem of Menger. This
notion also has a strong practical background. Let 
$G$ be a digraph of connectivity $k$ and diameter $D$.
By Menger's theorem, between any two distinct vertices
$x$ and $y$ in $G$ there are $k$ vertex-disjoint paths.
Such a set of disjoint paths, denoted by $C(x, y)$,
is called a container, and its length is defined
as maximum length among the paths in the container. 
The wide distance from $x$ to $y$ is then defined
to be the minimum length of the containers from $x$ to $y$,
and the wide diameter is the maximum wide distance among
all the pairs of distinct vertices. As we have mentioned
before, the consideration of the wide diameter of a network
has solid practical background which we will not
get into the discussions. Clearly,
the wide diameter of a digraph is at least as large
as the ordinary diameter. However, it is rather
remarkable that for most of the popular interconnection networks
the wide diameters are not significantly bigger
than (actually, a small constant bigger than)
the ordinary diameter, like the hypercube,
the de Bruijn, the Kautz,  and the star networks. 
The main result of this section is to show that
such a remarkable phenomenon also occurs in the
cycle prefix network. For a vertex $X$ in \cp,
we shall use $N(X)$ to denoted the set of
vertices next to $X$, and $M(X)$ the set of
vertices adjacent to $X$. We shall use the notation
$i\circ X$ to denote the vertex adjacent to $X$ that is
obtained by rotating the element $i$ to the beginning
position if $i$ is $X$, or by shifting $i$ into $X$
and bumping the last element out of $X$, namely $S_i(x)$
by the previous notation. If $i=x_1$, let $i\circ X = X$.
Note that $N(X)$ consists of
vertices $X_i = i\circ X$ for $i\not= x_1$, and $M(Y)$
consists of vertices 
\[ Y_i = \left\{ 
         \begin{array}{ll}
	    2\,3\,\cdots \, i \, 1 \, (i+1)\, \cdots \, D \, ,
			\quad &  \mbox{if}\quad 2\leq i \leq D\, . \\
          \equsep
	    2\, 3\, \cdots \, D \, i\, \quad  &
		\mbox{if}\quad D < i \leq \Delta+1\, ,\\
          \equsep Y&\mbox{if}\quad i = 1\,.

	    	  \end{array} \right.
\]

\begin{thm}
The wide diameter of \cp\ is at most $D+2$. It is exactly
$D+2$ for $D\geq 4$. Specifically, if $X$ and $Y$ are distinct
vertices in \cp, then there is a bijection $\theta$ of $N(X)\setminus{\{Y\}}$
to $M(Y)\setminus{\{X\}}$ such that the shortest paths from $Z$ to $\theta(Z)$
are vertex disjoint and do not contain either $X$ or $Y$.
\end{thm}

Since the proof of the above theorem heavily depends on
the unique shortest path property, we here give
a brief review of the shortest path routing in \cp.
Given two vertices $X$ and $Y$ in \cp, a tail of $Y$
with respect to $X$ (as the origin) is a suffix
$y_{k+1}\cdots y_D$ such that it forms a subsequence of $X$,
say $x_{i_1} x_{i_2} \cdots x_{i_{D-k}}$, and
all the elements $x_{1}, x_{2}, \ldots, x_{i_{D-k}}$
occur in $Y$. Note that a tail can be the empty sequence.
A header of $Y$ with respect to $X$ is a prefix
$y_1\cdots y_{k}$ such that the complement suffix $y_{k+1}\cdots y_D$
is a tail. It is proved in [\ref{FMC93}] that
the distance from $X$ to $Y$ is the length of the
shortest header of $Y$ with respect to $X$.
Suppose $y_1y_2\cdots y_{k}$ is the shortest header of $Y$
with respect to $X$, then the shortest path from $X$
to $Y$ is determined by the following prefixes:
\[ X \quad \Rightarrow \quad y_k **\* \quad \Rightarrow \quad
	  y_{k-1}y_{k} **\,* \quad \Rightarrow \quad \cdots \quad
	\Rightarrow \quad y_1\cdots y_k **** = Y\, ,\]
where $***$ is the usual wild-card notation to mean
``some sequence'' in order to fill the gap in the notation of
a sequence.  Without
loss of generality $Y$ can be assumed to be the 
standard origin.
For clarity we list the 
following conditions which together are equivalent to $d(X,Y)=k$ for
$k<D$:
\begin{itemize}
\item[(a).]   $y_D$ appears in $X$, say $x_j=y_D$.
\item[(b).]   $x_1, x_2, \ldots, x_j$ are all in $Y$ (but they
		do not necessarily form a subsequence).
\item[(c).]   $y_{k+1}\cdots y_D$ is a subsequence of $X$,
		but $y_k y_{k+1}$ is not. 
\end{itemize}
When $d(X, Y)=D$, it is equivalent to the following statement:
\begin{itemize}
\item[(a).]  either $y_D$ is not in $X$,
\item[(b).]  or $y_D$ is in $X$, say $x_j=y_D$, but there exists
	 	$x_r$ with $r<j$ that is not in $Y$.
\end{itemize}

For example, let $X=47285136$ and $Y=82164753$, the shortest header
of $Y$ with respect to $X$ is illustrated by $8216\,|\,4753$
The following property is
helpful in understanding the routing scheme in
\cp\ and it implies the uniqueness of the
shortest path.

\begin{prop}
Given $X=x_1x_2\cdots x_D$ and $Y=y_1y_2\cdots y_D$ 
in $\Gamma_\Delta(D)$, suppose
$k=d(X,Y)$. Then $d(i\circ X, Y) \geq d(X,Y)$ unless
$i=y_k$, in which case $d(i\circ X, Y)= d(X,Y)-1$.
\end{prop}

In order to reach the conclusion in the above
theorem concerning the wide diameter of \cp, we shall
start with the easiest case $x_1=1$. We
give a complete treatment of this case. For the other cases,
we only give an outline of the proof. The
details are similar to the case $x_1 = 1$ but more tedious.
In this regard, we hope that a simpler construction
will be found with a better understanding of the
the wide diameter of \cp. There is no doubt that 
the construction given in this
paper is {\em ad hoc}, although it does give
the best bound.

For the case $x_1=1$, the mapping $\theta$ is defined by
\[ \theta(X_i) = Y_i\,, \quad (2\leq i \leq \Delta+1).\]
The following is an example for $\Delta=5, D=4$ and $X=1325$.
\[ 
\begin{array}{lllllllllll}
X_2= & 2135   & \rightarrow & {\underline 4} 213  &
	 \rightarrow & \underline{34}21   &
	 \rightarrow & \underline{134} 2  &
	 \rightarrow & \underline{2134}  & = Y_2 \\
\equsep
X_3= & 3125   &	\rightarrow & \underline{4}312 &
		\rightarrow & \underline{14}32 &
		\rightarrow & \underline{314}2 &
		\rightarrow & \underline{2314} & = Y_3 \\
\equsep
X_4= & 4132   &	\rightarrow & \underline{3}412 &
		\rightarrow &  &
		            &  &
		            & \underline{23}41 &  = Y_4 \\
\equsep
X_5= & 5132   &	\rightarrow & \underline{4}513 &
		\rightarrow & \underline{34}51 &
		\rightarrow &
		    &     &  \underline{234}5  & = Y_5 \\
\equsep
X_6= & 6132   &	\rightarrow & \underline{4}613 &
		\rightarrow & \underline{34}61 &
		\rightarrow & 
 			 & 
		& \underline{234}6  & = Y_6 
\end{array}\]

The following lemma gives the distance from
$X_i$ to $Y_i$, from which the shortest path routing
is determined in terms of the shortest header.
 
\begin{lem}
Let $X=x_1x_2\cdots x_D$ be a vertex in \cp\ such that $x_1=1$. 
Suppose the distance from $X$ to $Y=12\cdots D$ is $k$,
then the distance from $X_i$ to $Y_i$  is given by
\[  d(X_i, Y_i) = \left\{
	\begin{array}{ll}
        k, \quad & \mbox{if} \quad 1 < i < k , \\
\equsep        k-2, \quad & \mbox{if} \quad  i= k , \\
\equsep        i-2, \quad & \mbox{if} \quad k < i \leq D , \\
\equsep        D-1, \quad & \mbox{if} \quad i > D\, .
\end{array}\right.
\]
\end{lem}

{\em Proof.} We first consider the case when $k<D$.
Note that $d(X, Y)=k$ is equivalent to the above conditions
(a), (b) and (c) altogether. We need to check the same
conditions for the corresponding distances in various cases.

For $1<i<k$, the verification for $d(i\circ X, Y_i)=k$ 
is divided into the following three steps:

(a). $Y_D=D$ appears in $X$: Suppose $i\circ X$ does not
contain $D$. Since $D$ is in $X$, it must be the last
element in $X$ and $i\circ X$ is obtained from $X$
via a shift operation. By the condition (b) for $d(X,Y)=k < D$,
$x_1, x_2, \ldots, x_D$ are all in $Y$, implying that
$X$ is a permutation on $1, 2, \ldots, D$. Thus, $i\circ X$
is obtained from $X$ via a rotation, which is a contradiction.

(b). Suppose $x_j=D$.  If $i \circ X$ is next to 
$X$ via a rotation, then every element prior to $D$ in $i\circ X$
is in $Y_i$ since $Y_i$ is a permutation of $Y$. If $i\circ X$
is next to $X$ via a shift operation, the above argument for (a)
shows that $D$ cannot be the last element of $X$. It also follows
that every element prior to $D$ in $i\circ X$ is still in $Y_i$.

(c). Since $(k+1, k+2, \ldots, D)$ is a subsequence of $X$,
it follows that it is also a subsequence of $i\circ X$ because
$i<k$ and $D$ stays in $i\circ X$. Clearly, $(k, k+1)$ cannot
be a subsequence of $i\circ X$ because it is not a subsequence of $X$.

Combining (a), (b) and (c) one sees that $d(X_i, Y_i)=k$, and the
shortest header of $Y_i$ with respect to $X_i$ is illustrated as 
follows:
\[ 2 \, 3 \cdots i \, 1 \,(i+1) \, \cdots k \, | \, (k+1) \cdots D\, .\]

For the case $i=k$, $d(X_k, Y_k)=k-2$: The verification of conditions
(a) and (b) are the same as for the previous case. The only catch
for condition (c) is that $(k, 1, k+1, \ldots D)$ is a subsequence
of $k\circ X$. Since $k$ is the first element of $k\circ X$, it follows
that $d(X_k,Y_k)=k-2$ and the shortest header of $Y_k$ is illustrated below:
\[ 2 \, 3\, \cdots \, (k-1) \, | \, k \, 1\, (k+1) \cdots D\, .\]

For the case $k< i \leq D$, $d(X_i, Y_i)=i-2$: The argument
for conditions (a) and (b) remain the same.  Noticing
that the first two elements of $i\circ X$ is $i\,1$ and that
$(i, 1, i+1, \cdots, D)$ is a subsequence of $i\circ X$, 
it follows that the shortest header of $Y_i$ is illustrated
below:
\[ 2\,3\, \cdots, k\, (k+1) \, | \, i \, 1\, (i+1)\, \cdots \, D\, .\]

It now comes the last subcase: $D< i \leq \Delta+1$. The
tail containing the single element $i$ of $Y_i$ is clearly
the longest tail with respect to $X_i$, and it is illustrated below:
\[ 2\, 3\, \cdots \, D \, | \, i \, .\]

We finally finish up the main case $k=D$, for which $Y$ is
not a closed vertex with respect to $X$. For $1<i<D$,
if $D$ is not in $X$, then $D$ is not in $i\circ X$ either.
Suppose $x_j=D$ and there exists $x_r$ with $r<j$ that is not in $Y$.
If $i\circ X$ does not contain $D$, then we are done. 
If $i\circ X$ contains $D$, then $x_r$ stays in $i\circ X$, but
it is not in $Y_i$. Hence we still have $d(X_i, Y_i)=D$.

For $i=D$, we have $D\circ X = D\, 1\, ***$ and $Y_D= 2\, 3\, \cdots \, D \,1$.
Clearly $d(X_D, Y_D)=D-2$.  For $i>D$, this is an easy matter, and the
same as for the case $k<D$. This completes all the cases. \qed

The shortest path routing from $X_i$ to $Y_i$ easily follows from the
above lemma. Our next goal is to show that all the shortest paths from 
$X_i$ to $Y_i$ are vertex-disjoint. To this end, we need to define
two statistics on a vertex so that they can be used to distinguish
the vertices along the shortest paths from $X_i$ to $Y_i$.
Given $X=x_1x_2\cdots x_D$, suppose $d(X, Y)=k$ where $Y=12\cdots D$.
Define $\alpha(X)$ to be the first element $x_i$ such that
$x_i\not\in \{ k+1, k+2,\ldots, D+1\}$. 
Let $\beta(X,i) = j+1$, where $j$ is obtained as follows:
Let $Y'$ be the second to the last vertex on the shortest path from
$X$ to $i\circ Y$. Let $j$ be the element immediately preceding
$i$ in $Y'$ or, if $j$ does not occur in $Y'$, the last element of $Y'$.
Equivalently, $\beta(X,i)$ is the minimum of $D+2$ and the
smallest $x > k$ such that $x$
is to the right of $i$ or not in $X$.
Let $\beta(X) = \beta(X,\alpha(X))$.
For example, suppose $X=531624$ and $Y=123456$.
Then $d(X,Y)=4$, $\alpha(X)=3$, $\beta(X)=6$, $\beta(X,2) = 7$.
We call
$(\alpha(x), \beta(x))$ the {\em characteristic pair} of $X$. 

It turns out that the characteristic pair of a vertex on
the shortest path from $X_i$ to $Y_i$ can be easily determined
along with the routing. The following table illustrates the
shortest path $P_i$ from $X_i$ to $Y_i$ in various cases,
together with the characteristic pairs from which one sees
that all the vertices are indeed distinct.
The notation $j\rightarrow$ means the operation of
getting $j\circ Z$ from $Z=z_1z_2\cdots z_D$ for $j\not= z_1$.

\begin{table}
\caption{ Case 1: $x_1 = 1$.}
\label{table:case1}%
\[\begin{array}{|llr|c|c|p{75pt}|}
\hline
i\hspace{20pt}\mbox{} & P_i&& \alpha(v) &  \beta(v) & Notes \\
\hline
\multicolumn{2}{|l}{\mbox{(a)} 1<i<k:}&&&&\\&
  \begin{array}{l}
    i\rightarrow\\
    k\rightarrow\\
    \vdots\\
    i+1\rightarrow\\
  \end{array}&\left.\begin{array}{r}\\\\\\\\\end{array}\right\}&
    i&  k+1 &  $i\not=1$  \\
  &\begin{array}{l}
    1\rightarrow\\
    i\rightarrow\\
    \vdots\\
    2\rightarrow\\
  \end{array}&\left.\begin{array}{r}\\\\\\\\\end{array}\right\} &
    1 & i+1 & $i+1\not=k+1$ \\
\hline
\multicolumn{2}{|l}{\mbox{(b) } i=k:}&&&&\\&
  \begin{array}{l}
    k \rightarrow\\
    \vdots\\
    2\rightarrow\\
  \end{array}&\left.\begin{array}{r}\\\\\\\end{array}\right\} &
  1 & k+1 & 
\\
\hline
\multicolumn{2}{|l}{\mbox{(c)} k+1\leq i\leq D+1:}&&&&\\&
   \begin{array}{l}
     i\rightarrow\\
     \vdots\\
     2\rightarrow\\
   \end{array}&\left.\begin{array}{r}\\\\\\\end{array}\right\}&
  1 & i+1 & \\
\hline
\multicolumn{2}{|l}{\mbox{(d)} D+1<i\leq \Delta+1:}&&&&\\&
   \begin{array}{l}
     i\rightarrow\\
     D\rightarrow\\
     \vdots\\
     2\rightarrow\\
   \end{array}&\left.\begin{array}{r}\\\\\\\\\end{array}\right\}&
    i & D+1 & \\
\hline
\end{array}
\]
\end{table}

The cases other than $x_1=1$ are more tedious.
The critical part is to construct
the mapping from $N(x)$ to $M(Y)$. Recall that $d(X,Y)=k$.

For $x_1=k+1$, we have
\[ \theta(i\circ X) = \left\{ 
     \begin{array}{ll}
	       Y_k, \quad   & \mbox{if} \quad i=1\, ,\\
\equsep   	Y_i, \quad  & \mbox{if} \quad 1< i < k\; , \\
\equsep 	Y_{\beta(X,1)-1} \quad & \mbox{if} \quad i=k\, ,\\
\equsep  	Y_{i-1}  \quad   & \mbox{if} \quad k+1\leq i<\beta(X,1)\, ,\\
\equsep    	Y_i \quad & \mbox{if} \quad \beta(X,1)\leq i \leq \Delta+1\,.
\end{array} \right. \]

It is straightforward to see that $\theta$ is a bijection. 
The distance from $i\circ X$ to $\theta(i\circ X)$ are given
below:
\[      d(i\circ X, \theta(i \circ X)) = \left\{
	\begin{array}{ll}
	        k-1, \quad   & \mbox{if} \quad i=1\, ,\\
\equsep   	k-1, \quad  & \mbox{if} \quad 1< i < k\; , \\
\equsep 	k-1, \quad & \mbox{if} \quad i=k\, ,\\
\equsep  	i,  \quad   & \mbox{if} \quad k+1\leq i<\beta(X,1)\, ,\\
\equsep    	D-1 \quad & \mbox{if} \quad \beta(X,1)\leq i \leq \Delta+1\,.
\end{array} \right. \]
We omit the detailed verification of the above distances. Based
on these distances, we have the following table which illustrates
the shortest path routing from $i\circ X$ to $\theta(i\circ X)$. 
In addition to the characteristic pairs, we need one more
characteristic to distinguish the vertices. For the sake of
easy description, we assume that the elements in $X$ that 
are greater than $D$ occur in increasing order starting
$D+1, D+2, \ldots$, because a permutation on the set 
$\{ D+1, D+2, \ldots, \Delta+1\}$ can map the vertex $X$ into this
form without affecting the destination vertex or the other
elements in $X$ that do not exceed $D$.

\begin{table}
\caption{ Case 2: $x_1 = k+1$.}
\label{table:case2}%
\[\begin{array}{|llr|c|c|c|p{75pt}|}
\hline
i\hspace{20pt}\mbox{} & P_i&& \alpha(v) & \beta(v) &\beta(v, 1)& Notes\\
\hline
\multicolumn{2}{|l}{\mbox{(a)} i = 1:}&&&&&\\&
  \begin{array}{l}
    1\rightarrow\\
    k \rightarrow\\
    \vdots\\
    2\rightarrow\\
   \end{array}&\left.\begin{array}{r}\\\\\\\\\end{array}\right\}&
    1& k+1 & k+1 & Ignore if $k+1 = 2$.\\
\hline
\multicolumn{2}{|l}{\mbox{(b)} 1<i<k:}&&&&&\\&
  \begin{array}{l}
    i \rightarrow\\
    k \rightarrow\\
    \vdots\\
    i+1\rightarrow\\
  \end{array}&\left.\begin{array}{r}\\\\\\\\\end{array}\right\}&
    i& k  & \beta(x,1)&\\&
  \begin{array}{l}
    1\rightarrow\\
    i\rightarrow\\
    \vdots\\
    2\rightarrow\\
  \end{array}&\left.\begin{array}{r}\\\\\\\\\end{array}\right\}&
    1& i+1& i+1& $i+1\not=k+1$\\
\hline
\multicolumn{2}{|l}{\mbox{(c)} i = k:}&&&&&\\&
  \begin{array}{l}
    k\rightarrow\\
    \vdots\\
    2\rightarrow\\
  \end{array}&\left.\begin{array}{r}\\\\\\\end{array}\right\}&
    \not=k+1 & k+1&\beta(x, 1)&Ignore if $k+1 = 2$.\\
\hline
\multicolumn{2}{|l}{\mbox{(d)} k+1< i < \beta(x,1):} &&&&&\\&
   \begin{array}{l}
     i\rightarrow\\
   \end{array}&\left.\begin{array}{r}\\\end{array}\right\}&
     k+1 & i+1 & \beta(x,1) & \\ &
   \begin{array}{l}
     1\rightarrow\\
     i-1\rightarrow\\
     \vdots\\
     2\rightarrow\\
   \end{array}&\left.\begin{array}{r}\\\\\\\\\end{array}\right\}&
     1&i&i&\\
\hline
\multicolumn{2}{|l}{\mbox{(e)} \beta(x,1)\leq i\leq D+1:}&&&&&\\&
   \begin{array}{l}
     i\rightarrow\\
     \vdots\\
     2\rightarrow\\
   \end{array}&\left.\begin{array}{r}\\\\\\\end{array}\right\}&
     &i+1&i+1&$i+1\not=\beta(x, 1)$\\
\hline
\multicolumn{2}{|l}{\mbox{(f)} D+1<i\leq \Delta+1, i\not=k+1:}&&&&&\\&
   \begin{array}{l}
     i\rightarrow\\
     D\rightarrow\\
     \vdots\\
     2\rightarrow\\
   \end{array}&\left.\begin{array}{r}\\\\\\\\\end{array}\right\}&
     i&D+1&D+1&\\
\hline
  \end{array}\]
\end{table}

For $x_1\not= 1$ or $k+1$, we have
\[ \theta(i\circ X) = \left\{ 
        \begin{array}{ll}
	    Y_{\beta(X,1)-1}, \quad   & \mbox{if} \quad i=1\, ,\\
\equsep   	Y_i, \quad  & \mbox{if} \quad 1< i < k\; , \\
\equsep 	Y_{\alpha(X)} \quad & \mbox{if} \quad i=k\, ,\\
\equsep  	Y_{i-1}  \quad   & \mbox{if} \quad k+1\leq i<\beta(X,1)\, ,\\
\equsep    	Y_i \quad & \mbox{if} \quad \beta(X,1)\leq i \leq \Delta+1\,.
\end{array} \right. \]
In this case, the distances are given below:
\[      d(i\circ X, \theta(i \circ X)) = \left\{
	\begin{array}{ll}
	        \beta(X,1)-1, \quad   & \mbox{if} \quad i=1\, ,\\
\equsep   	k-1, \quad  & \mbox{if} \quad 1< i < k\; , i\not= x_1, \\
\equsep 	k-2, \quad & \mbox{if} \quad i=k\, ,\\
\equsep  	i,  \quad   & \mbox{if} \quad k+1\leq i<\beta(X,1)\, ,\\
\equsep         i-2, \quad  & \mbox{if} \quad \beta(X,1)\leq i \leq D\,,\\
\equsep    	D-1 \quad & \mbox{if} \quad D < i \leq \Delta+1, i\not=k+1\,.
\end{array} \right. \]
Note that the assumption on the elements in $X$ that are greater than
$D$ implies that $x_1<k$.  
The detailed information on the shortest path routing from $i\circ X$
to $\theta(i\circ X)$ is given in the table below, which also includes
the statistic $\beta(V,1)$ to distinguish the vertices.

\begin{table}
\caption{ Case 3.}
\label{table:case3}%
\[\begin{array}{|llr|c|c|c|p{75pt}|}
\hline
\hspace{20pt}\mbox{} & P_i&& \alpha(v) & 
	\beta(v) & \beta(v,1)& Notes\\
\hline
\multicolumn{2}{|l}{\mbox{(a)} i=1:}&&&&&\\&
  \begin{array}{l}
      1\rightarrow\\
  \end{array}&\left.\begin{array}{r}\\\end{array}\right\}&
      1 & k+1 & k+1 & \\&
  \begin{array}{l}
      \beta(x,1)-1\rightarrow\\
      \vdots\\
      2 \rightarrow\\
  \end{array}&\left.\begin{array}{r}\\\\\\\end{array}\right\}&
      1 & \beta(x,1) & \beta(x,1)&\\
\hline
\multicolumn{2}{|l}{\mbox{(b) } 1 < i < k, i\not=x_1:}&&&&&\\&
   \begin{array}{l}
      i\rightarrow\\
      k\rightarrow\\
      \vdots\\
      i+1\rightarrow\\
    \end{array}&\left.\begin{array}{r}\\\\\\\\\end{array}\right\}&
      i & k+1 &\beta(x,1)&$i\not=x_1$\\&
    \begin{array}{l}
      1\rightarrow\\
      i\rightarrow\\
      \vdots\\
      2\rightarrow\\
    \end{array}&\left.\begin{array}{r}\\\\\\\\\end{array}\right\}&
      1& i+1&i+1&$i+1\not=k+1$, $i+1\not = x_1+1$\\
\hline
\multicolumn{2}{|l}{\mbox{(c) } i = k:}&&&&&\\&
   \begin{array}{l}
      k \rightarrow\\
      \vdots\\
      x_1 + 1\rightarrow\\
   \end{array}&\left.\begin{array}{r}\\\\\\\end{array}\right\}&
     x_1& k+1 &\beta(x,1)& $d(v,Y)<k$\\&
   \begin{array}{l}
      1\rightarrow\\
      x_1\rightarrow\\
      \vdots\\
      2\rightarrow\\
   \end{array}&\left.\begin{array}{r}\\\\\\\\\end{array}\right\}&
     1& x_1+1& x_1+1 &$x_1<k$\\
\hline
\multicolumn{2}{|l}{\mbox{(d) } k+1 \leq i <\beta(x,1):}&&&&&\\&
    \begin{array}{l}
      i\rightarrow\\
    \end{array}&\left.\begin{array}{r}\\\end{array}\right\}&
      x_1& i+1&\beta(x,1)&$i+1\not=k+1$\\&
    \begin{array}{l}
      1\rightarrow\\
      i-1\rightarrow\\
      \vdots\\
      2\rightarrow\\
    \end{array}&\left.\begin{array}{r}\\\\\\\\\end{array}\right\}&
      1&i&i&\\
\hline
\multicolumn{2}{|l}{\mbox{(e) } \beta(x,1)\leq i\leq D+1:}&&&&&\\&
    \begin{array}{l}
      i\rightarrow\\
      \vdots\\
      2\rightarrow\\
    \end{array}&\left.\begin{array}{r}\\\\\\\end{array}\right\}&
      x_1&i+1&i+1&\\
\hline
\multicolumn{2}{|l}{\mbox{(f) }D+1<i\leq \Delta+1, i\not=k+1:}&&&&&\\&
    \begin{array}{l}
     i\rightarrow\\
     D\rightarrow\\
     \vdots\\
     2\rightarrow\\
    \end{array}&\left.\begin{array}{r}\\\\\\\\\end{array}\right\}&
     i&D+1&D+1&\\
\hline
\end{array}\]
\end{table}

After completing the above case by case analysis, we arrive at the
conclusion that the wide diameter of $\Gamma_\Delta(D)$ does not
exceed $D+2$. To determine when the wide diameter is exactly $D+2$,
let $X=(\Delta+1, \Delta, \ldots, \Delta-D+2)$ and $Y=12\cdots D$.
The distance from $i\circ Y$ is $D$ except for $i=D$.
If $1< i < D$, the unique shortest path from $i\circ X$ to $Y$
goes through $23\cdots D \, (\Delta+1)$. Provided that $D-2\geq 2$,
this implies that, of any $\Delta$ disjoint path from $X$ to $Y$,
at least one must have length at least $D+2$.  \qed

We have left open the problem of finding the
wide diameter of \cpr, which is of great interest
if it is determined. 

\vspace{1cm}

\noindent
{\small {\large \bf Acknowledgments.} 
This work was performed under the auspices of the U. S. Department of Energy.
We thank D. F. Hsu for helpful discussions.
}

\vspace{1cm}

\newpage

\noindent
{\Large\bf References}

\newcounter{paper}
\setcounter{paper}{0}
\begin{list}
{[\arabic{paper}]\hfill}{\usecounter{paper} 
   \setlength{\leftmargin=1.5cm}{\labelsep=2mm}}

\item \label{AK89} 
	S. B. Akers and B. Krishnamurthy,
	A group-theoretic model for symmetric interconnection
	networks, IEEE Trans. on Computers, 38 (1989), 
	555-566.

\item 	\label{Ber92} 
	J.-C. Bermond, ed., Interconnection Networks, 
	Special Issue in Discrete Applied Mathematics, Vol. 37/38 (1992).

\item \label{CFK93}
       W. Y. C. Chen, V. Faber and E. Knill,
       A new routing scheme for cycle prefix graphs,
       LAUR-93-3576, Los Alamos National Laboratory,
	1993.

\item \label{CG}
	J. H. Conway and M. J. T. Guy,
	Message graphs, Ann. of Discrete Math., 13 (1982), 
	61-64.

\item \label{FMC93}
       V. Faber, J. Moore and W. Y. C. Chen,	
	Cycle prefix digraphs for symmetric interconnection networks,
	Networks, 23 (1993), 641-649.

\item   \label{CF92}  F. Comellas and M. A. Fiol,
	Vertex symmetric digraphs with small diameter, preprint,  1992.

\item \label{Hsu93}
	D. F. Hsu, ed., Interconnection Networks and Algorithms,
	Special Issue in Networks, to appear.

\item \label{Hsu**}
	D. F. Hsu, On container width and length in graphs,
	groups, and networks, to appear.

\item \label{JR92}
        M. Jiang and F. Ruskey,
	Determining the Hamilton-Connectedness of certain
	vertex-transitive graphs, Technical Report, DCS-202-IR,
	Department of Computer Science, University of Victoria,
        Victoria, B.C., Canada, 1992.

\item 	\label{Sab69}
	G. Sabidussi,
	Vertex transitive graphs,
	Monatsh. Math.
	68 (1969), 426-438.

\end{list}

\end{document}